\documentclass{amsart}
\usepackage{amsfonts}
\usepackage{graphicx}
\usepackage{amscd}

\newtheorem{theorem}{Theorem}[section]
\theoremstyle{plain}

\newtheorem{lemma}[theorem]{Lemma}

\newtheorem{proposition}[theorem]{Proposition}

\theoremstyle{example}
\newtheorem{example}[theorem]{Example}
\theoremstyle{remark}
\newtheorem{remark}[theorem]{Remark}
\numberwithin{equation}{section}

\input{tcilatex}

\begin{document}
\title[Minimal hypersurfaces]{On Deformable Minimal Hypersurfaces in Space
Forms}
\author{Andreas Savas-Halilaj}
\address{Department of Mathematics and Statistics, University of Cyprus,
P.O. Box 20537, 1678 Nicosia, Cyprus}
\email{halilaj.andreas@ucy.ac.cy}
\date{}

\begin{abstract}
The aim of this paper is to complete the local classification of minimal
hypersurfaces with vanishing Gauss-Kronecker curvature in a 4-dimensional
space form. Moreover, we give a classification of complete minimal
hypersurfaces with vanishing Gauss-Kronecker curvature and scalar curvature
bounded from below.
\end{abstract}

\maketitle

\section{Introduction}

A classical result stated by R. Beez \cite{Beez}, but completely proved by
W. Killing \cite{Killing}, says that a hypersurface in the Euclidean space $%
\mathbb{R}^{n+1}$ is locally rigid if its type number, i.e., the rank of its
Gauss map, is at least $3$. Hence, the interesting case is that of
hypersurfaces with type number less than $3$.

Inspired by works of Sbrana \cite{Sbrana} and Cartan \cite{Cartan}, M.
Dajczer and D. Gromoll \cite{Dajczer-Gromoll} proved that every hypersurface
with constant type number in $\mathbb{R}^{n+1}$ can be described in terms of
its Gauss map. This kind of description is called the \textquotedblleft
Gauss parametrization\textquotedblright . However, the Gauss parametrization
provides a description of the hypersurface at the points where the type
number is constant and there is no discussion on the manner in which
hypersurfaces with different type numbers can be produced or joined together.

In this paper we consider minimal hypersurfaces with type number less than $%
3 $ in $\mathbb{Q}_{c}^{4}$. Our aim is to provide a local parametrization
of such hypersurfaces around the totally geodesic points. In the $3$%
-dimensional case the assumption that the type number is less than $3$ is
equivalent to the assumption that the Gauss-Kronecker curvature vanishes. It
turns out that hypersurfaces with constant type number $2$ are highly
deformable since they allow locally an one-parameter family of isometric
deformations, the so called associated family. We shall show that such
hypersurfaces are either totally geodesic or arise, at least locally, as
normal bundles of certain branched minimal surfaces. More precisely, we
shall show the following result.

\begin{theorem}
Let $M^{3}$ be an oriented Riemannian manifold of dimension $3$ and $%
f:M^{3}\rightarrow \mathbb{Q}_{c}^{4}$ a minimal isometric immersion with
vanishing Gauss-Kronecker curvature. Then, $f\left( M^{3}\right) $ is either
totally geodesic in $\mathbb{Q}_{c}^{4}$ or, at least locally, coincide with
the polar map of a branched\textit{\ minimal surface.}
\end{theorem}

\noindent For the definition of the polar map over branched minimal surfaces
see Section 4.

Dajczer and Gromoll \cite{Dajczer-Gromoll} showed that in order to solve the
classification problem for the class of complete minimal hypersurfaces with
type number less than $3$, one has to understand the nature of $3$%
-dimensional complete minimal hypersurfaces with vanishing Gauss-Kronecker
curvature in $\mathbb{Q}_{c}^{4}$.

In \cite{Almeida} S.C. de Almeida and F.G.B. Brito initiated the study of
compact minimal hypersurfaces in $\mathbb{S}^{4}$ with vanishing
Gauss-Kronecker curvature and gave a classification under the assumption
that the second fundamental form of the hypersurface nowhere vanishes. J.
Ramanathan in \cite{Ramanathan} removed the assumption on the second
fundamental form and proved that such hypersurfaces arise as normal bundles
over branched minimal surfaces in $\mathbb{S}^{4}$, unless they are totally
geodesic. We point out, that Theorem 1.1 generalizes Ramanathan's result,
since the same conclusion holds without the assumption of compactness.

Complete minimal hypersurfaces with zero Gauss-Kronecker curvature in $%
\mathbb{Q}_{c}^{4}$ was investigated by T. Hasanis et al. in \cite{Hasanis1,
Hasanis2, Hasanis3}. In particular, a classification of such hypersurfaces
was given, under the assumptions that the scalar curvature is bounded from
below and that the second fundamental form of $f$ is nowhere zero. In the
following theorem, we generalize results in \cite{Hasanis1, Hasanis2,
Hasanis3, Ramanathan}.

\begin{theorem}
\textit{Let }$M^{3}$ \textit{be an oriented, }$3$\textit{-dimensional,
complete Riemannian manifold with scalar curvature }bounded from below. If $%
f:M^{3}\rightarrow \mathbb{Q}_{c}^{4}$ is\ \textit{a minimal isometric
immersion with Gauss-Kronecker curvature identically zero, then either }$f$
is totally geodesic,\textit{\ or}

\begin{enumerate}
\item[$a)$] $c=0$ and\textit{\ }$f\left( M^{3}\right) $ \textit{splits as a
Euclidean product }$M^{2}\times \mathbb{R}$\textit{, where }$M^{2}$\textit{\
is a complete minimal surface in }$\mathbb{R}^{3}$,

\item[$b)$] $c=1$ and $f\left( M^{3}\right) $ \textit{is the image of the
polar map associated with a branched superminimal immersion} $%
g:M^{2}\rightarrow \mathbb{S}^{4}$,

\item[$c)$] $c=-1$ and $f\left( M^{3}\right) $ is the hyperbolic cylinder
over \textit{a complete minimal surface }$M^{2}$ in $Q^{3}$\textit{, where }$%
Q^{3}$ is\textit{\ a horosphere or\ an equidistant hypersurface in} $\mathbb{%
H}^{4}$.
\end{enumerate}
\end{theorem}

The paper is organized as follows. In Section 2 we introduce the necessary
background about harmonic maps. In Section 3 we review some basic facts
about the behavior of the Gauss map of branched minimal surfaces. In Section
4 we furnish a method to produce minimal hypersurfaces in $\mathbb{Q}%
_{c}^{4} $ with zero Gauss-Kronecker curvature. Section 5 is devoted to the
study of the set of totally geodesic points of such hypersurfaces. Finally,
in Section 6 we the prove Theorem 1.1 and Theorem 1.2.

\section{Harmonic Maps}

Let $M^{m}$ and $N^{n}$ be two differentiable manifolds, of dimensions $m$
and $n$, equipped with Riemannian metrics $\left\langle \text{ },\text{ }%
\right\rangle _{M}$ and $\left\langle \text{ },\text{ }\right\rangle _{N}$,
respectively. Denote by $\nabla ^{M}$ and $\nabla ^{N}$ the corresponding
Levi-Civita connections. If $F:M^{m}\rightarrow N^{n}$ is a smooth map, we
define the \textit{energy integral} of $F$ over a compact domain $\Omega $
of $M^{m}$ by%
\begin{equation*}
E\left( F,\Omega \right) :=\frac{1}{2}\int_{\Omega }\left\Vert dF\right\Vert
^{2}dM,
\end{equation*}%
where $dM$ is the volume element of $M^{m}$ and $\left\Vert dF\right\Vert
\left( p\right) $ is the Hilbert-Schmidt norm of the differential $%
dF:TM^{m}\rightarrow TN^{n}$ of $F$ at a point $p$ of $\Omega $. A smooth
map $F$ is said to be \textit{harmonic} if it is a critical point of the
energy integral $E\left( \cdot ,\Omega \right) $, with respect to variations
supported in $\Omega $.

The Hessian $\nabla ^{2}F:TM^{m}\times TM^{m}\rightarrow F^{\ast }\left(
TN^{n}\right) $ of $F$ is defined by%
\begin{equation*}
\nabla ^{2}F\left( X,Y\right) :=\nabla _{X}^{F}dF\left( Y\right) -dF\left(
\nabla _{X}^{M}Y\right) ,\quad X,Y\in TM^{m},
\end{equation*}%
where $F^{\ast }\left( TN^{n}\right) $ stands for the induced vector bundle
by $F$ and $\nabla ^{F}$ for the induced Levi-Civita connection of $F^{\ast
}\left( TN^{n}\right) $. Taking the Euler-Lagrange equations, it follows
that $F$ is harmonic if and only if the trace with respect to $\left\langle 
\text{ },\text{ }\right\rangle _{M}$ of the \textit{Hessian} $\nabla ^{2}F$
vanishes identically. In local coordinates $\left( U,\phi \right) $ on $%
M^{m} $ and $\left( V,\psi \right) $ on $N^{n}$, the harmonicity of $F$ is
expressed by the following system of elliptic differential equations%
\begin{equation}
\Delta f^{\alpha }+\dsum\limits_{\alpha ,\beta =1}^{n}\Gamma _{\beta \gamma
}^{\alpha }\circ F\left\langle \nabla ^{M}f^{\beta },\nabla ^{M}f^{\gamma
}\right\rangle _{M}=0,\quad \alpha \in \left\{ 1,2,\cdots ,n\right\} ,
\end{equation}%
where $\Delta $ is the Laplace operator on $M^{m}$, $f^{\alpha }$ the
components of $F$ and $\Gamma _{\beta \gamma }^{\alpha }$ the Christoffel
symbols of $N^{n}$ with respect to the coordinate system $\left( V,\psi
\right) $.

The notion of harmonic map can be carried out successfully even in the case
of semi Riemannian manifolds without any changes. However, in this article
we shall need only the case of harmonic maps $F:M^{2}\rightarrow N_{1}^{n}$,
where $N_{1}^{n}$ is a semi Riemannian manifold of signature $\left(
1,n-1\right) $. We point that in this case the system of differential
equations $\left( 2.1\right) $ is still elliptic.

An isometric immersion between Riemannian manifolds is harmonic if and only
if it is minimal. Moreover, any harmonic map between real analytic manifolds
with real analytic Riemannian metrics is also real analytic (cf., \cite[p.
117]{Eells-Sampson}). According to the so called Unique Continuation
Property, states that two harmonic maps $F_{1},F_{2}:M^{m}\rightarrow N^{n}$
which\ coincide in an open subset $\Omega $ of $M^{n}$ should coincide in
all of $M^{m}$ (cf., \cite[Theorem 1, p. 213]{Sampson}).

Let assume that $\mathcal{K}$ is a compact set in some open ball $%
B_{R}\left( 0\right) \subset \mathbb{R}^{n}$ with radius $R$ centered at the
origin of $\mathbb{R}^{n}$. Then, the \textit{relative }$2$\textit{-capacity}
of $\mathcal{K}$ is defined to be%
\begin{equation*}
\limfunc{cap}\nolimits_{2,R}\left( \mathcal{K}\right) :=\inf \left\{ \int_{%
\mathbb{R}^{n}}\left\vert \nabla \varphi \right\vert ^{2}d\mathbb{R}%
^{n}:\varphi \in C_{c}^{\infty }\left( B_{R}\left( 0\right) \right) \text{
and }\varphi \left( x\right) \geq 1,\text{ for all }x\in \mathcal{K}\right\}
,
\end{equation*}%
where $C_{c}^{\infty }\left( B_{R}\left( 0\right) \right) $ is the space of
smooth functions $\varphi :B_{R}\left( 0\right) \rightarrow \mathbb{R}$ with
compact support. The relative 2-capacity measures how \textit{%
\textquotedblleft small\textquotedblright\ }must be a set $\mathcal{K}$ in
order a harmonic map to be extendable over it. For $n\geq 3$, we set $%
\limfunc{cap}_{2}\left( \mathcal{K}\right) =\limfunc{cap}_{2,\infty }\left( 
\mathcal{K}\right) $. If $\mathcal{C}$ is an arbitrary subset of $\mathbb{R}%
^{n}$ we say that $\limfunc{cap}_{2}\left( \mathcal{C}\right) =0$ if and
only if $\limfunc{cap}_{2}\left( \mathcal{K}\right) =0$, for any compact
subset $\mathcal{K}$ of $\mathcal{C}$. In particular, if $M^{m}$ is an $m$%
-dimensional differentiable manifold and $\mathcal{C}$ a subset of $M^{m}$,
then we say that $\limfunc{cap}_{2}\left( \mathcal{C}\right) =0$ if and only
if for any coordinate chart $\left( U,\phi \right) $ on $M^{m}$ it holds $%
\limfunc{cap}_{2}\left( \phi \left( \mathcal{C}\cap U\right) \right) =0$.
Finally, there is a close relation between the $2$-capacity and the $\left(
m-2\right) $-Hausdorff measure $\mathcal{H}^{m-2}$. In fact, the following
result is well known (cf., \cite[pp. 154-158]{Evans-Gariepy}).

\begin{theorem}
Let $\mathcal{K}$ be a compact set. Then, the following assertions are valid:

\begin{enumerate}
\item[$a)$] If $\mathcal{H}^{m-2}\left( \mathcal{K}\right) <\infty $, then $%
\limfunc{cap}_{2}\left( \mathcal{K}\right) =0$.

\item[$b)$] If $\limfunc{cap}_{2}\left( \mathcal{K}\right) =0$, then $%
\mathcal{H}^{s}\left( \mathcal{K}\right) =0$ for any $s>m-2$.
\end{enumerate}
\end{theorem}

\noindent For example, if $\mathcal{C}$ is a picewise smooth curve in $%
\mathbb{R}^{3}$, then $\mathcal{H}^{1}\left( \mathcal{C}\right) =1$ and
consequently $\limfunc{cap}_{2}\left( \mathcal{C}\right) =0$. M. Meier in 
\cite{Meier} provides a criterion about the removability of singularities
for harmonic maps between Riemannian manifolds.

\begin{theorem}
Let $M^{m}$ be a $m$-dimensional Riemannian manifold, $\mathcal{C}$ a
relatively closed subset of $M^{m}$ with $\limfunc{cap}_{2}\left( \mathcal{C}%
\right) =0$ and $N^{n}$ a $n$-dimensional Riemannian manifold. If $%
F:M^{m}\rightarrow N^{n}$ is a map satisfying the following properties:

\begin{enumerate}
\item[$a)$] $F$ is continuous on $M^{m}$,

\item[$b)$] $F$ is harmonic on $M^{m}-\mathcal{C}$.
\end{enumerate}

\noindent then, $F$ extends analytically to a harmonic map in all of $M^{m}$.
\end{theorem}

\section{The Gauss map of branched surfaces}

In this section we review some basic facts about the Gauss map of spacelike
branched minimal surfaces in semi Riemannian manifolds.

Denote by $\mathbb{R}_{s}^{n}$ the a $n$-dimensional vector space $\mathbb{R}%
^{n}$ endowed with the semi Riemannian metric $\left\langle \text{ },\text{ }%
\right\rangle _{s}$ given by%
\begin{equation*}
\left\langle v,w\right\rangle _{s}:=\left( -1\right)
^{s}v_{1}w_{1}+v_{2}w_{2}+\cdots +v_{n}w_{n},\quad s\in \left\{ 0,1\right\} ,
\end{equation*}%
where $v=\left( v_{1},v_{2},\cdots ,v_{n}\right) $ and $w=\left(
w_{1},w_{2},\cdots ,w_{n}\right) $. The \textit{Grassmannian space }of the $%
2 $-dimensional planes in $\mathbb{R}_{s}^{n}$ is denoted by $\mathbb{G}%
_{2}\left( \mathbb{R}_{s}^{n}\right) $ and admits a differentiable structure
of dimension $2\left( n-2\right) $. A 2-plane $\Pi $ in $\mathbb{R}_{s}^{n}$
is said to be \textit{spacelike}, if the restriction to $\Pi $ of the inner
product $\left\langle \text{ },\text{ }\right\rangle _{s}$ is positive
definite. Let $\mathbb{G}_{2}^{+}\left( \mathbb{R}_{s}^{n}\right) $ the
space of all oriented spacelike 2-planes in $\mathbb{R}_{s}^{n}$. Consider a 
$2$-dimensional oriented Riemannian manifold $\left( M^{2},\left\langle 
\text{ },\text{ }\right\rangle \right) $ endowed with a fixed complex
structure as a Riemann surface, a $n$-dimensional semi Riemannian manifold $%
\left( N_{s}^{n},\left\langle \text{ },\text{ }\right\rangle _{s}\right) $
equipped with a semi Riemannian metric of signature $\left( s,n-s\right) $, $%
s\in \left\{ 0,1\right\} $. Assume that $N_{s}^{n}$ can be isometrically
immersed into a Minkowski space $\mathbb{R}_{s}^{n+l}$. The latest
assumption is not essential for $s=0$ since, according to a well known
result of J. Nash, any Riemannian manifold can be isometrically immersed
into a Euclidean space.

We say that $g:M^{2}\rightarrow N_{s}^{n}\subset \mathbb{R}_{s}^{n+l}$ is a 
\textit{(weakly)\ conformal map}, if%
\begin{equation*}
g^{\ast }\left\langle \text{ },\text{ }\right\rangle _{s}=E\left\langle 
\text{ },\text{ }\right\rangle ,
\end{equation*}%
where $E:M^{2}\rightarrow \mathbb{R}$ is a smooth non negative function, the
so called \textit{conformal factor}. On the open set $M_{\ast }^{2}$ of its
regular points, $g$ is a spacelike immersion and its Gauss map $\xi :M_{\ast
}^{2}\rightarrow \mathbb{G}_{2}^{+}\left( \mathbb{R}_{s}^{n+l}\right) $ is
defined by%
\begin{equation*}
\xi \left( x\right) :=dg\left( T_{x}M_{\ast }^{2}\right) ,\quad x\in M^{2}.
\end{equation*}%
If $x_{0}$ is a critical point of $g$, then the tangent plane of $g$ at $%
x_{0}$ does not exists in general. A critical point $x_{0}$ of $g$ is called 
\textit{branch point of order }$m$, $m\in \mathbb{N}$, if there are
coordinate charts $\left( U,\varphi \right) $ centered at $x_{0}$ and $%
\left( V,\psi \right) $ around $g\left( x_{0}\right) $, such that the
representation map $\widetilde{g}:\varphi \left( U\right) \mathbb{%
\rightarrow C\times R}^{n-2}$, $\widetilde{g}:=\psi \circ g\circ \varphi
^{-1}$, has the form%
\begin{equation*}
\widetilde{g}\left( z\right) =\left( z^{m+1}+h_{1}\left( z\right) +\sqrt{-1}%
h_{3}\left( z\right) ,h\left( z\right) \right) ,\quad z\in \varphi \left(
U\right) ,
\end{equation*}%
where $h_{1}\left( z\right) $, $h_{2}\left( z\right) $, $h\left( z\right)
\in O\left( \left\vert z\right\vert ^{m+2}\right) $. In particular, $g$ is
called \textit{branched surface} if each critical point is a branch point.
According to results due to R.D. Gulliver et al. \cite[Lemma 1.3 and Lemma
3.1]{Osserman}, the branch points are isolated and the Gauss map can be
continuously extended over the branch points. It turns out that the Gauss
map of a conformal harmonic map can be smoothly extended over the branch
points. More precisely, the following theorem due to J.-H. Eschenburg and R.
Tribuzy \cite{Eschenburg-Tribuzy} provides a characterization of the
critical points of a conformal harmonic map.

\begin{theorem}
Let $M^{2}$ be a Riemann surface and $N_{s}^{n}$ a semi Riemannian manifold
of signature $\left( s,n-s\right) $, $s\in \left\{ 0,1\right\} $, which can
be isometrically immersed into a Minkowski space $\mathbb{R}_{s}^{n+l}$.
Assume that $g:M^{2}\rightarrow N_{s}^{n}\subset \mathbb{R}_{s}^{n+l}$ is a
conformal harmonic map. Then,

\begin{enumerate}
\item[$a)$] the critical points of $g$ are branch points.

\item[$b)$] The Gauss map of $g$ extends smoothly over the branch points.

\item[$c)$] If $x_{0}\in M^{2}$ is a branch point of order $m\geq 1$, then
for any holomorphic chart $\left( U,z\right) $, $z:U\rightarrow \mathbb{C}$,
centered at $x_{0}$, we have%
\begin{equation*}
E\left( z\right) =\left\vert z\right\vert ^{2m}E_{0}\left( z\right) ,
\end{equation*}%
where $E_{0}$ is a smooth positive function on $U$.
\end{enumerate}
\end{theorem}

\begin{remark}
Although the proof in \cite{Eschenburg-Tribuzy} was given for surfaces lying
in a Riemannian manifold, the same technique can be applied also in the case
of spacelike surfaces in a semi Riemannian manifold $N_{1}^{n}$ of signature 
$\left( 1,n-1\right) $, which can be immersed isometrically in $\mathbb{R}%
_{1}^{n+l}$.
\end{remark}

\section{Gauss parametrization}

In this section we shall describe a method of constructing minimal
hypersurfaces with vanishing Gauss-Kronecker curvature in a space form $%
\mathbb{Q}_{c}^{4}$ and will illustrate the results with some examples.

Consider a branched minimal surface $g:M^{2}\rightarrow \mathbb{S}%
^{3}\subset \mathbb{R}^{4}$, where $M^{2}$ is an oriented manifold of
dimension $2$. Denote by $\left\langle \text{ },\text{ }\right\rangle _{g}$
the induced by $g$ metric tensor on $M^{2}$. On the open set $M_{\ast }^{2}$
of its regular points, $g$ is an immersion and $\left\langle \text{ },\text{ 
}\right\rangle _{g}$ is a Riemannian metric. According to Theorem 3.1, the
normal bundle%
\begin{equation*}
\mathcal{N}\left( g\right) :=\left\{ \left( x,w\right) \in M^{2}\times 
\mathbb{R}^{4}:w\bot \mathbb{R\cdot }g+dg\left( T_{x}M^{2}\right) \right\} ,
\end{equation*}%
is a well defined smooth manifold of dimension 3. Let $\gamma
:M^{2}\rightarrow \mathbb{R}$ be a smooth function\ and denote by $\nabla
\gamma $ its gradient, by $\nabla ^{2}\gamma $ its Hessian and by $\Delta
\gamma $ its Laplacian with respect to $\left\langle \text{ },\text{ }%
\right\rangle _{g}$ on $M_{\ast }^{2}$. Assume further that $dg\left( \nabla
\gamma \right) $ can be smoothly extended over the branch points of $g$.
Then, the map $\Psi _{\gamma ,g}:\mathcal{N}\left( g\right) \rightarrow 
\mathbb{R}^{4}$ given by%
\begin{equation*}
\Psi _{\gamma ,g}\left( x,w\right) :=\gamma \left( x\right) g\left( x\right)
+dg\left( \nabla \gamma \right) \left( x\right) +w,
\end{equation*}%
is well defined and smooth. The map $\Psi _{\gamma ,g}$, when restricted on
the normal bundle of $g:M_{\ast }^{2}\rightarrow \mathbb{S}^{3}$, coincides
with the \textit{polar map} invented by Dajczer and Gromoll in \cite[Section
1]{Dajczer-Gromoll}. Since our considerations are of local nature, we shall
replace $M^{2}$ with the open unit disc%
\begin{equation*}
B:=\left\{ z\in \mathbb{C}:\left\vert z\right\vert <1\right\} .
\end{equation*}%
From now on denote by $B_{\ast }:=B-\left\{ 0\right\} $ the punctured open
unit disc.

\begin{proposition}
Let $g:B\rightarrow \mathbb{S}^{3}\subset \mathbb{R}^{4}$ be a branched
minimal surface, having $z=0$ as branch point of order $m$ and being
everywhere else\ regular. Assume that $\gamma :B\rightarrow \mathbb{R}$ is a
smooth function, satisfying on $B_{\ast }$\ the differential equation%
\begin{equation*}
\Delta \gamma +2\gamma =0.
\end{equation*}%
Then,

\begin{enumerate}
\item[$a)$] the polar map is well defined and smooth if and only if $%
dg\left( \nabla \gamma \right) $ can be smoothly extended over the branch
point.

\item[$b)$] Let $A_{w}$ be the shape operator of $g$ with respect to the
normal $w$. The map $\Psi _{\gamma ,g}$ is regular at $\left( z,w\right) $, $%
z\neq 0$, if and only if the self adjoint operator $P:T_{z}B\rightarrow
T_{z}B$ given by%
\begin{equation*}
P:=\nabla ^{2}\gamma +\gamma I-A_{w},
\end{equation*}%
is non-singular. Moreover, $\Psi _{\gamma ,g}$ is regular at $\left(
0,w\right) $ if and only if there is a negative number $l$ such that%
\begin{equation*}
\lim_{z\rightarrow 0}\left\vert z\right\vert ^{2m}\det P=l.
\end{equation*}

\item[$c)$] On the open set of its regular points, the polar map is minimal
with vanishing Gauss-Kronecker curvature.
\end{enumerate}
\end{proposition}

\begin{proof}
Since $z=0$ is a branch point, according to Theorem 3.1, there is a positive
smooth function $E_{0}$ on $B$ such that%
\begin{equation*}
\left\langle \text{ },\text{ }\right\rangle _{g}=\left\vert z\right\vert
^{2m}E_{0}\left( dx^{2}+dy^{2}\right) .
\end{equation*}%
Denote by $\eta $ a unit normal vector field along $g$. Parametrize $%
\mathcal{N}\left( g\right) $ by $B\times \mathbb{R}$ via the map $\left(
z,t\right) \rightarrow \left( z,t\eta \left( z\right) \right) $. Then, for
any $\left( z,t\right) \in B_{\ast }\times \mathbb{R}$,%
\begin{equation*}
\Psi _{\gamma ,g}\left( z,t\right) =\gamma \left( z\right) g\left( z\right)
+dg\left( \nabla \gamma \right) \left( z\right) +t\eta \left( z\right) .
\end{equation*}%
$a)$ It is obvious that the polar map is well defined and smooth if and only
if $dg\left( \nabla \gamma \right) $ can be smoothly extended over the
branch point.

\noindent $b)$ On an arbitrary point $\left( z,t\right) $, $z\neq 0$, we have%
\begin{equation*}
d\Psi _{\gamma ,g}\left( \partial _{t}\right) =\eta
\end{equation*}%
and, for any $X\in T_{z}B$,%
\begin{equation}
d\Psi _{\gamma ,g}\left( X\right) =dg\left( PX\right) +\omega \left(
X\right) \eta ,
\end{equation}%
where $\omega $ stands for the dual form of the vector $A_{\eta }\nabla
\gamma $. The metric tensor $\left\langle \text{ },\text{ }\right\rangle
_{\Psi }$ induced on $B_{\ast }\times \mathbb{R}$ by $\Psi _{\gamma ,g}$ is
given by the formula%
\begin{equation*}
\left\langle \text{ },\text{ }\right\rangle _{\Psi }=\left\langle
P^{2},\right\rangle _{g}+\left( dt+\omega \right) \otimes \left( dt+\omega
\right) .
\end{equation*}%
Since $\gamma $ satisfies the differential equation $\Delta \gamma +2\gamma
=0$, we deduce that%
\begin{eqnarray*}
\left\langle \text{ },\text{ }\right\rangle _{\Psi } &=&-\det P\left\langle 
\text{ },\text{ }\right\rangle _{g}+\left( dt+\omega \right) \otimes \left(
dt+\omega \right) \\
&=&-\left\vert z\right\vert ^{2m}E_{0}\det P\left( dx^{2}+dy^{2}\right)
+\left( dt+\omega \right) \otimes \left( dt+\omega \right) .
\end{eqnarray*}%
Hence, the polar map is regular at $\left( z,t\right) $, $z\neq 0$, if and
only if $\det P\neq 0$. Since,%
\begin{equation*}
\omega \left( X\right) =\left\langle A_{\eta }\nabla \gamma ,X\right\rangle
_{g}=-\left\langle d\eta \left( X\right) ,dg\left( \nabla \gamma \right)
\right\rangle ,
\end{equation*}%
it follows that $\omega $ is well defined and smooth over the branch point.
Moreover, $\Psi _{\gamma ,g}$ is regular at $\left( 0,t\right) $ if and only
if there exists a negative real number $l$, such that%
\begin{equation*}
\lim_{z\rightarrow 0}\left\vert z\right\vert ^{2m}\det P=l.
\end{equation*}%
$c)$ The unit vector field $\xi $ given by $\xi \left( z,t\right) :=g\left(
z\right) $ is smooth and normal along $\Psi _{\gamma ,g}$. Suppose now that $%
\left( 0,t_{0}\right) $ is a regular point of $\Psi _{\gamma ,g}$. Then,
there exists a neighborhood $U$ around $\left( 0,t_{0}\right) $ where $\Psi
_{\gamma ,g}$ is regular. Since,%
\begin{equation*}
\lim_{z\rightarrow 0}\left\vert z\right\vert ^{2m}\det P=l<0,
\end{equation*}%
there exists a continuous positive function $\Phi :U\rightarrow \mathbb{R}$,
such that%
\begin{equation*}
-\left\vert z\right\vert ^{2m}\det P=\Phi \left( z,t\right) ,\quad \left(
z,t\right) \in U.
\end{equation*}%
By $\left( 4.1\right) $, it follows that the principal curvatures of $\Psi
_{\gamma ,g}$ on $U$ are given by%
\begin{equation*}
k_{1}\left( z,t\right) =-k_{3}\left( z,t\right) =\frac{\left\vert
z\right\vert ^{m}}{\sqrt{\Phi \left( z,t\right) }},\text{ \ }k_{2}\left(
z,t\right) =0.
\end{equation*}%
and the proof of the theorem is complete.
\end{proof}

\begin{example}
Let $G:B\rightarrow \mathbb{S}^{3}\subset \mathbb{R}^{4}$ be a minimal
immersion, $\left\langle \text{ },\text{ }\right\rangle _{G}$ the induced
Riemannian metric on $B$ by $G$ and $A_{G}$ the shape operator of $G$ with
respect to a unit normal vector field $\eta $ along $G$.\ The map $%
g:B\rightarrow \mathbb{S}^{3}$ which assigns to each point of $B$ the
endpoint of $\eta $, when this is translated to the origin of $\mathbb{R}%
^{4} $, coincide with the Gauss map of $G$ and is harmonic. Furthermore, the
induced metric tensor $\left\langle \text{ },\text{ }\right\rangle $ by $g$
on $B$ is given by the formula%
\begin{equation*}
\left\langle \text{ },\text{ }\right\rangle =-\det A_{G}\left\langle \text{ }%
,\text{ }\right\rangle _{G}.
\end{equation*}%
According to \cite[ Lemma 1.4, p. 338]{Lawson}, the zeroes of $\det A$ are
isolated and of even order, unless $G$ is totally geodesic and $\det A$
vanishes identically. Suppose that $z=0$ is an isolated totally geodesic
point of $G$. Then, there exists a smooth positive function $\Lambda $ such
that%
\begin{equation*}
-\det A_{G}=\left\vert z\right\vert ^{2m}\Lambda ^{2},\quad z\in B.
\end{equation*}%
Thus, $z=0$ is a branch point of $g$ with order $m$. Denote by $\nabla $ the
Levi-Civita connection of $\left\langle \text{ },\text{ }\right\rangle $ and
by $A$ the shape operator of $g$ on $B_{\ast }$. From the Weingarten formula
it follows that $A_{G}A=I$. Let $\alpha $ be a constant vector of $\mathbb{R}%
^{4}$ and $\gamma ,\sigma :B\rightarrow \mathbb{R}$ the functions given by
the formulas $\gamma :=\left\langle g,\alpha \right\rangle $ and $\sigma
:=\left\langle G,\alpha \right\rangle $, respectively. Then, for any $z\in
B_{\ast }$,%
\begin{equation*}
dg\left( \nabla \gamma \right) =\alpha ^{G}\text{ \ and \ }\nabla ^{2}\gamma
=\sigma A-\gamma I,
\end{equation*}%
where $\alpha ^{G}$ is the tangential component of $\alpha $ along $G$.
Hence, $dg\left( \nabla \gamma \right) $ is well defined and smooth over the
branch point $z=0$. Moreover, $\gamma $ satisfies the differential equation $%
\Delta \gamma +2\gamma =0$ on $B_{\ast }$, where $\Delta $ stands for the
Laplacian with respect to $\left\langle \text{ },\text{ }\right\rangle $.
Consider now the polar map $\Psi _{\gamma ,g}:B\times \mathbb{R\rightarrow R}%
^{4}$,%
\begin{equation*}
\Psi _{\gamma ,g}\left( z,t\right) =\gamma \left( z\right) g\left( z\right)
+dg\left( \nabla \gamma \right) \left( z\right) +tG\left( z\right) .
\end{equation*}%
For any $z\in B_{\ast }$, we have%
\begin{equation*}
\det P=\frac{\left( \sigma -t\right) ^{2}}{\det A_{G}}=-\frac{\left( \sigma
-t\right) ^{2}}{\left\vert z\right\vert ^{m}\Lambda ^{2}}.
\end{equation*}%
Consequently, the polar map fails to be regular at points $\left( z,t\right) 
$ with $t=\sigma \left( z\right) $. On the open set of its regular points, $%
\Psi _{\gamma ,g}$ has principal curvatures%
\begin{equation*}
k_{1}\left( z,t\right) =-k_{3}\left( z,t\right) =\frac{\left\vert
z\right\vert ^{m}\Lambda \left( z\right) }{\left\vert \sigma \left( z\right)
-t\right\vert },\quad k_{2}\left( z,t\right) =0.
\end{equation*}%
Notice that points of the form $\left( 0,t\right) $, $t\neq \sigma \left(
z\right) $, corresponds to totally geodesic points of the polar map.
\end{example}

\begin{example}
Assume that $g:B\subset \mathbb{S}^{2}\rightarrow \mathbb{S}^{3}$ is totally
geodesic immersion and $\gamma :B\rightarrow \mathbb{R}$ a smooth function
satisfying $\Delta \gamma +2\gamma =0$. Then, the polar map $\Psi _{\gamma
,g}$ is a cylinder over a minimal surface lying in $\mathbb{R}^{3}$.
Conversely, the cylinder over a complete minimal surface in $\mathbb{R}^{3}$%
, gives rise to a complete minimal hypersurface with vanishing
Gauss-Kronecker curvature in $\mathbb{R}^{4}$. The cylinders are the only
known examples of complete minimal hypersurfaces with Gauss-Kronecker
identically zero in $\mathbb{R}^{4}$.
\end{example}

Assume that $M^{2}$ is an oriented manifold of dimension $2$ and $%
g:M^{2}\rightarrow \mathbb{S}^{4}\subset \mathbb{R}^{5}$ a branched minimal
surface. Then, the unit normal bundle%
\begin{equation*}
\mathcal{N}_{1}\left( g\right) :=\left\{ \left( x,w\right) \in M^{2}\times 
\mathbb{R}^{5}:\left\vert w\right\vert =1,\text{ }w\bot \mathbb{R\cdot }%
g+dg\left( T_{x}M^{2}\right) \right\}
\end{equation*}%
is a well defined smooth manifold of dimension $3$. The map $\Psi _{g}:%
\mathcal{N}_{1}\left( g\right) \rightarrow \mathbb{S}^{4}$ given by%
\begin{equation*}
\Psi _{g}\left( x,w\right) :=w,
\end{equation*}%
is called the \textit{polar map} over $g$.

Following the notation of Proposition 4.1, we obtain the following result.

\begin{proposition}
Let $g:B\rightarrow \mathbb{S}^{4}\subset \mathbb{R}^{5}$ be a branched
minimal surface, having $z=0$ as a branch point of order $m$ and being
everywhere else regular. Then,

\begin{enumerate}
\item[$a)$] a point $\left( z,w\right) $, $z\neq 0$, is a regular point of $%
\Psi _{g}$ if and only if the shape operator $A_{w}$ of $g$, with respect to
the normal $w$, is non-singular. Moreover, $\Psi _{g}$ is regular at a point
of the form $\left( 0,w\right) $, if and only if there exists a negative
number $l$ such that%
\begin{equation*}
\lim_{z\rightarrow 0}\left\vert z\right\vert ^{2m}\det A_{w}=l.
\end{equation*}

\item[$b)$] On the open set of its regular points, the polar map is minimal
with vanishing Gauss-Kronecker curvature.
\end{enumerate}
\end{proposition}

\begin{proof}
Since $z=0$ is a branch point of $g$, the metric tensor $\left\langle \text{ 
},\text{ }\right\rangle _{g}$ can be expressed as%
\begin{equation*}
\left\langle \text{ },\text{ }\right\rangle _{g}=\left\vert z\right\vert
^{2m}E_{0}\left( dx^{2}+dy^{2}\right) ,
\end{equation*}%
where $E_{0}$ is a positive smooth function on $B$. Denote by $\eta _{3}$, $%
\eta _{4}$ an orthonormal frame field along $g$ and parametrize $\mathcal{N}%
_{1}\left( g\right) $ by $B\times \mathbb{S}^{1}$ via the map $\left(
z,t\right) \rightarrow \left( z,\cos t\eta _{3}\left( z\right) +\sin t\eta
_{4}\left( z\right) \right) $. Then, $\Psi _{g}$ can be represented as%
\begin{equation*}
\Psi _{g}\left( z,t\right) =\cos t\eta _{3}\left( z\right) +\sin t\eta
_{4}\left( z\right) .
\end{equation*}%
Notice that the connection form%
\begin{equation*}
\omega _{34}\left( X\right) :=\left\langle d\eta _{3}\left( X\right) ,\eta
_{4}\right\rangle
\end{equation*}%
is well defined and smooth even at the branch point of $g$.

\noindent $a)$ Let $\left( z,t\right) $, $z\neq 0$, be a point of $B\times 
\mathbb{S}^{1}$. Then,%
\begin{equation*}
d\Psi _{g}\left( \partial _{t}\right) =-\sin t\eta _{3}+\cos t\eta _{4},
\end{equation*}%
and, for any $X\in T_{z}B$,%
\begin{equation}
d\Psi _{g}\left( X\right) =-dg\left( A_{w}X\right) +\omega _{34}\left(
X\right) d\Psi _{g}\left( \partial _{t}\right) .
\end{equation}%
The induced by $\Psi _{g}$ metric tensor $\left\langle \text{ },\text{ }%
\right\rangle _{\Psi }$ on $B_{\ast }\times \mathbb{S}^{1}$ is given by the
formula%
\begin{equation*}
\left\langle \text{ },\text{ }\right\rangle _{\Psi }=-\left\vert
z\right\vert ^{2m}E_{0}\det A_{w}\left( dx^{2}+dy^{2}\right) +\left(
dt+\omega _{34}\right) \otimes \left( dt+\omega _{34}\right) .
\end{equation*}%
Hence, $\left( z,t\right) $, $z\neq 0$, is a regular point of $\Psi _{g}$ if
and only if the shape operator $A_{w}$ of $g$ with respect to the normal $w$
is non-singular. Moreover, $\left( 0,t\right) $ is a regular point of $\Psi
_{g}$ if and only if there is a negative real number $l$ such that%
\begin{equation*}
\lim_{z\rightarrow 0}\left\vert z\right\vert ^{2m}\det A_{w}=l.
\end{equation*}%
$b)$ The unit vector field $\xi \left( z,t\right) :=g\left( z\right) $ is
smooth and normal along $\Psi _{g}$. Assume that $\left( 0,t_{0}\right) $ is
a regular point of the polar map. Then, there exists a neighborhood $U$
around $\left( 0,t_{0}\right) $, where $\Psi _{g}$ is regular. Since,%
\begin{equation*}
\lim_{z\rightarrow 0}\left\vert z\right\vert ^{2m}\det A_{w}=l<0,
\end{equation*}%
there exists a continuous positive function $\Phi :U\rightarrow \mathbb{R}$,
such that%
\begin{equation*}
-\left\vert z\right\vert ^{2m}\det A_{w}=\Phi \left( z,t\right) ,\quad
\left( z,t\right) \in U.
\end{equation*}%
Then, by $\left( 4.2\right) $, we deduce that $\Psi _{g}$ has principal
curvatures%
\begin{equation*}
k_{1}\left( z,t\right) =-k_{3}\left( z,t\right) =\frac{\left\vert
z\right\vert ^{m}}{\sqrt{\Phi \left( z,t\right) }},\text{ \ }k_{2}\left(
z,t\right) =0,
\end{equation*}%
and the proof is complete.
\end{proof}

\begin{example}
Let $G:B\rightarrow \mathbb{S}^{3}\subset \mathbb{S}^{4}$ be\ a minimal
immersion having $z=0$ as an isolated\ totally geodesic point of order $2m$.
Then, its Gauss map $g:M^{2}\rightarrow \mathbb{S}^{3}$ is a branched
minimal surface, having $z=0$ branch point of order $m$. If $A_{G}$ is the
shape operator of $G$, then there exists a positive smooth function $\Lambda 
$ such that%
\begin{equation*}
-\det A_{G}=\left\vert z\right\vert ^{2m}\Lambda ^{2}\text{.}
\end{equation*}%
Consider $\mathbb{S}^{3}$ as a totally geodesic hypersurface of $\mathbb{S}%
^{4}$ and let $\alpha $ be\ the unit normal vector of $\mathbb{S}^{3}$ in $%
\mathbb{S}^{4}$. The polar map $\Psi _{g}:B\times \mathbb{S}^{1}\rightarrow 
\mathbb{S}^{4}$ has the form%
\begin{equation*}
\Psi _{g}\left( z,t\right) =\cos t\alpha +\sin tG\left( z\right) .
\end{equation*}%
The determinant of the shape operator $A_{w}$ of $g$, with respect to a unit
normal $w=\cos t\alpha +\sin tG$, is equal to%
\begin{equation*}
\det A_{w}=\frac{\cos ^{2}t}{\det A_{G}}=-\frac{\cos ^{2}t}{\left\vert
z\right\vert ^{2m}\Lambda ^{2}}.
\end{equation*}%
Thus, $\Psi _{g}$ fails to be an immersion at points $\left( z,t\right) $
with $t=\left( 2n+1\right) \pi /2$, $n\in \mathbb{Z}$. Furthermore, on the
open set of its regular points, $\Psi _{g}$ has principal curvatures%
\begin{equation*}
k_{1}\left( z,t\right) =-k_{3}\left( z,t\right) =\frac{\left\vert
z\right\vert ^{m}\Lambda \left( z\right) }{\left\vert \cos t\right\vert },%
\text{ \ }k_{2}\left( z,t\right) =0.
\end{equation*}%
Notice that points of the form $\left( 0,t\right) $, $t\neq \left(
2n+1\right) \pi /2$, $n\in \mathbb{Z}$, corresponds to totally geodesic
points of the polar map.
\end{example}

\begin{example}
R.L. Bryant studied in \cite{Bryant} an interesting class of branched
minimal surfaces in $\mathbb{S}^{4}$, the so called superminimal ones.
Actually, a branched minimal surface is called superminimal if at each point
the ellipse of curvature is a circle or a point. Moreover, Bryant has shown
that all superminimal surfaces in $\mathbb{S}^{4}$ admit a \textquotedblleft
Weierstrass type representation\textquotedblright\ in terms of meromorphic
functions. More precisely, let $\mathbb{CP}^{3}$ be the\ complex projective
space endowed with the Fubini-Study metric and $\mathcal{P}:\mathbb{CP}%
^{3}\rightarrow \mathbb{S}^{4}$ the Penrose twistor fibration. Consider two
meromorphic functions $\phi $ and $\psi $ defined on a Riemannian surface $%
M^{2}$ and let $G:M^{2}\rightarrow \mathbb{CP}^{3}$ be the map given by%
\begin{equation*}
G\left( z\right) :=\left[ 1:\phi -\frac{1}{2}\psi \frac{d\phi }{d\psi }:\psi
:\frac{1}{2}\frac{d\phi }{d\psi }\right] .
\end{equation*}%
Then $G$ is holomorphic, horizontal with respect to the Penrose fibration
and $\mathcal{P}\circ G$ gives rise to a superminimal surface in $\mathbb{S}%
^{4}$. Conversely, any superminimal surface in $\mathbb{S}^{4}$ can be
produced, at least locally, in this way. Applying the polar map over
superminimal surfaces we can produce examples of minimal hypersurfaces (even
compact ones) in $\mathbb{S}^{4}$ with vanishing Gauss-Kronecker curvature.
For example, consider the pair of meromorphic functions $\phi ,\psi :\mathbb{%
CP}^{1}=\mathbb{C\cup }\left\{ \infty \right\} \rightarrow \mathbb{C}$,
given by $\phi \left( z\right) =z^{5}$ and $\psi \left( z\right) =z^{2}$.
Composing the holomorphic curve $G:\mathbb{CP}^{1}\rightarrow \mathbb{CP}%
^{3} $,%
\begin{equation*}
G\left( z\right) =\left[ 1:-\frac{1}{4}z^{5}:z^{2}:\frac{5}{4}z^{3}\right] ,
\end{equation*}%
with the Penrose fibration, we obtain a superminimal surface having at $z=0$
and $z=\infty $ branch points of order $1$. Moreover, it turns out that the
polar map over $\mathcal{P}\circ G$ is everywhere regular and the points of
the fibers over $z_{1}=0$ and $z_{2}=\infty $ corresponds to totally
geodesic points of the polar map.
\end{example}

Consider the hyperbolic space $\mathbb{H}^{4}$ as the hyperquadric%
\begin{equation*}
\mathbb{H}^{4}:=\left\{ x=\left( x_{1},x_{2},x_{3},x_{4},x_{5}\right) \in 
\mathbb{R}_{1}^{5}:\left\langle x,x\right\rangle _{1}=-1,\quad
x_{1}>0\right\} .
\end{equation*}%
The hyperquadric%
\begin{equation*}
\mathbb{S}_{1}^{4}:=\left\{ x\in \mathbb{R}_{1}^{5}:\left\langle
x,x\right\rangle _{1}=1\right\} ,
\end{equation*}%
is called the \textit{de Sitter space}. Let $M^{2}$ be an oriented manifold
of dimension 2 and $g:M^{2}\rightarrow \mathbb{S}_{1}^{4}\subset \mathbb{R}%
_{1}^{5}$ a branched spacelike stationary immersion. The timelike unit
normal bundle $\mathcal{N}_{-1}\left( g\right) $ of $g$, defined by%
\begin{equation*}
\mathcal{N}_{-1}\left( g\right) :=\left\{ \left( x,w\right) \in M^{2}\times 
\mathbb{R}_{1}^{5}:\left\langle w,w\right\rangle _{1}=-1,\text{ }w\bot 
\mathbb{R\cdot }g+dg\left( T_{x}M^{2}\right) \right\}
\end{equation*}%
is a smooth manifold of dimension $3$ and the \textit{polar map} $\Psi _{g}:%
\mathcal{N}_{-1}\left( g\right) \rightarrow \mathbb{H}^{4}$ is defined by
the formula%
\begin{equation*}
\Psi _{g}\left( x,w\right) :=w.
\end{equation*}%
Proceeding in the same way as in Proposition 4.4, we can prove the following
result.

\begin{proposition}
Let $g:B\rightarrow \mathbb{S}_{1}^{4}\subset \mathbb{R}_{1}^{5}$ be a
branched spacelike stationary surface, having $z=0$ as a branch point of
order $m$ and being everywhere else regular. Then,

\begin{enumerate}
\item[$a)$] a point $\left( z,w\right) $, $z\neq 0$, is a regular point of $%
\Psi _{g}$ if and only if the shape operator $A_{w}$ of $g$, with respect to
the normal $w$, is non-singular. Moreover, $\Psi _{g}$ is regular at a point
of the form $\left( 0,w\right) $, if and only if there exists a negative
number $l$ such that%
\begin{equation*}
\lim_{z\rightarrow 0}\left\vert z\right\vert ^{2m}\det A_{w}=l<0.
\end{equation*}

\item[$b)$] On the open set of its regular points, the polar map is minimal
with vanishing Gauss-Kronecker curvature.
\end{enumerate}
\end{proposition}

\begin{example}
Surfaces in $\mathbb{S}^{3}$ which are critical points of the functional%
\begin{equation*}
\mathcal{W}\left( M^{2}\right) =\int_{M^{2}}\left( H^{2}+1-K\right) dM,
\end{equation*}%
where $H$, $K$ are the mean and Gaussian curvature of the surface $M^{2}$,
are called Willmore surfaces. Bryant showed in \cite{Bryant1} that there is
a duality between branched spacelike stationary surfaces in $\mathbb{S}%
_{1}^{4}$ and Willmore surfaces in $\mathbb{S}^{3}$. More precisely, let $%
h:M^{2}\rightarrow \mathbb{S}^{3}\subset \mathbb{S}_{1}^{4}$ be an isometric
immersion of a $2$-dimensional oriented Riemannian manifold $M^{2}$.
Consider $\mathbb{S}^{3}$ as the submanifold of $\mathbb{S}_{1}^{4}$%
\begin{equation*}
\mathbb{S}^{3}=\left\{ \left( 0,y\right) \in \mathbb{R}_{1}^{5}=\mathbb{R}%
_{1}^{1}\times \mathbb{R}^{4}:\left\langle y,y\right\rangle _{1}=1\right\} ,
\end{equation*}%
and define the conformal Gauss map $\mathcal{G}_{h}:M^{2}\rightarrow \mathbb{%
S}_{1}^{4}$,%
\begin{equation*}
\mathcal{G}_{h}\left( x\right) :=H\left( x\right) \left( 1,h\left( x\right)
\right) +\left( 0,\eta \left( x\right) \right) ,\quad x\in M^{2},
\end{equation*}%
where $\eta $ is\ a unit normal vector field along $h$. The induced by $%
\mathcal{G}_{h}$ metric tensor $\left\langle \text{ },\text{ }\right\rangle
_{\sim }$ on $M^{2}$ is related with the Riemannian metric $\left\langle 
\text{ },\text{ }\right\rangle $ of $M^{2}$ by the formula%
\begin{equation*}
\left\langle \text{ },\text{ }\right\rangle _{\sim }=\left( H^{2}+1-K\right)
\left\langle \text{ },\text{ }\right\rangle ,
\end{equation*}%
where $K$ stands for the Gaussian curvature and $H$ for the mean curvature
of $h$. We point out that the zeroes of $H^{2}+1-K$ are isolated, unless $g$
is totally umbilical. Moreover, $\mathcal{G}_{h}$ is a branched spacelike
stationary surface if and only if $h$ is a Willmore surface. In fact,
branched spacelike stationary surfaces in $\mathbb{S}_{1}^{4}$ arise as
conformal Gauss maps of a Willmore surfaces in $\mathbb{S}^{3}$. Minimal
surfaces in $\mathbb{S}^{3}$ are the simplest examples of Willmore surfaces.
\end{example}

\begin{example}
Examples of complete minimal hypersurfaces in $\mathbb{H}^{4}$ arise as
hyperbolic cylinders over complete minimal surfaces in horospheres or
equidistant hypersurfaces. More precisely, let $Q^{3}$ be a horosphere or
equidistant hypersurface in $\mathbb{H}^{4}$ and denote by $\eta $ a unit
normal vector field of $Q^{3}$ in $\mathbb{H}^{4}$. If $h:M^{2}\rightarrow
Q^{3}\subset \mathbb{H}^{4}$ is a complete minimal surface, then the
hyperbolic cylinder $F_{h}:M^{2}\times \mathbb{R\rightarrow H}^{4}$, given by%
\begin{equation*}
F_{h}\left( z,t\right) :=\cosh th\left( z\right) +\sinh t\eta \left(
z\right) ,\ \ \left( z,t\right) \in M^{2}\times \mathbb{R},
\end{equation*}%
provides a complete minimal hypersurface in $\mathbb{H}^{4}$ with vanishing
Gauss-Kronecker curvature. In fact, these are the only known examples of
complete minimal hypersurface in $\mathbb{H}^{4}$ with vanishing
Gauss-Kronecker curvature.
\end{example}

\section{The nullity distribution}

Let $M^{3}$ be a $3$-dimensional, oriented Riemannian manifold and $%
f:M^{3}\rightarrow \mathbb{Q}_{c}^{4}$ an isometric immersion. Denote by $%
\xi $ a unit normal vector field along $f$ and by $A$ its shape operator.
From the Weingarten formula it follows that%
\begin{equation*}
d\xi =-df\circ A.
\end{equation*}%
The eigenvalues%
\begin{equation*}
k_{1}\geq k_{2}\geq k_{3}
\end{equation*}%
of $A$ are called \textit{principal curvatures} of the immersion. Recall
that the \textit{mean curvature} $H$, the \textit{square }$S$\textit{\ of
the second fundamental form} and the \textit{Gauss-Kronecker} curvature $K$
of the immersion are given in terms of the principal curvatures by the
following formulas%
\begin{equation*}
H=\frac{1}{3}\left( k_{1}+k_{2}+k_{3}\right) ,\quad
S=k_{1}^{2}+k_{2}^{2}+k_{3}^{2},\quad K=k_{1}k_{2}k_{3}.
\end{equation*}

Let us assume that $f$ is minimal and that the type number of the immersion
is at most $2$. Equivalently, $f$ is minimal with Gauss-Kronecker curvature
identically zero. A point $x_{0}$ where the square $S$ of the second
fundamental form vanishes is called a \textit{totally geodesic point.}

Assume at first that $x$ is a non-totally geodesic point of $f$. Then, there
is a simply connected neighborhood $\Omega $ around $x$, where the second
fundamental form of $f$ is nowhere zero. Thus, the principal curvatures of $%
f $ on $\Omega $ are $k_{1}=\lambda $, $k_{2}=0$, $k_{3}=-\lambda $, where $%
\lambda $ is a smooth positive function on $\Omega $. Moreover, we can
choose on $\Omega $ an orthonormal frame field $\left\{
e_{1},e_{2},e_{3}\right\} $ of principal directions corresponding to $%
\lambda $, $0$, $-\lambda $, respectively. Let $\left\{ \omega _{1},\omega
_{2},\omega _{3}\right\} $ and $\left\{ \omega _{ij}\right\} $, $i,j\in
\left\{ 1,2,3\right\} $, be the dual coframe and the connection forms,
respectively. Define the functions $u,v:U\rightarrow \mathbb{R}$, given by
the formulas%
\begin{equation*}
u=\left\langle \nabla _{e_{3}}e_{1},e_{2}\right\rangle \text{ \ and \ }%
v=\left\langle \nabla _{e_{1}}e_{1},e_{2}\right\rangle .
\end{equation*}%
From the Gauss and Codazzi equations, follows that%
\begin{equation}
\begin{array}{lll}
\omega _{12}\left( e_{1}\right) =v, & \omega _{13}\left( e_{1}\right) =\frac{%
1}{2}e_{3}\left( \log \lambda \right) \medskip , & \omega _{23}\left(
e_{1}\right) =u, \\ 
\omega _{12}\left( e_{2}\right) =0, & \omega _{13}\left( e_{2}\right) =\frac{%
1}{2}u,\medskip & \omega _{23}\left( e_{2}\right) =0, \\ 
\omega _{12}\left( e_{3}\right) =u, & \omega _{13}\left( e_{3}\right) =-%
\frac{1}{2}e_{1}\left( \log \lambda \right) , & \omega _{23}\left(
e_{3}\right) =-v,%
\end{array}%
\end{equation}%
and%
\begin{equation}
\begin{array}{ll}
e_{2}\left( v\right) =v^{2}-u^{2}+c, & e_{1}\left( u\right) =e_{3}\left(
v\right) , \\ 
e_{2}\left( u\right) =2uv, & e_{3}\left( u\right) =-e_{1}\left( v\right) .%
\end{array}%
\end{equation}%
The above equations yield%
\begin{eqnarray}
\left[ e_{1},e_{2}\right] &=&-ve_{1}+\tfrac{1}{2}ue_{3},\ \left[ e_{2},e_{3}%
\right] =\tfrac{1}{2}ue_{1}+ve_{3},\medskip \\
\left[ e_{1},e_{3}\right] &=&-\tfrac{1}{2}e_{3}\left( \log \lambda \right)
e_{1}-2ue_{2}+\tfrac{1}{2}e_{1}\left( \log \lambda \right) e_{3}.  \notag
\end{eqnarray}%
The following observation is a fact of basic importance.

\begin{lemma}
The functions $u$ and $v$ are harmonic.
\end{lemma}

\begin{proof}
From the equations $\left( 5.1\right) $ and the definition of the Laplacian,
we obtain%
\begin{eqnarray*}
\Delta u &=&e_{1}e_{1}\left( u\right) +e_{2}e_{2}\left( u\right)
+e_{3}e_{3}\left( u\right) \\
&&-\nabla _{e_{1}}e_{1}\left( u\right) -\nabla _{e_{2}}e_{2}\left( u\right)
-\nabla _{e_{3}}e_{3}\left( u\right) \\
&=&e_{1}e_{1}\left( u\right) +e_{2}e_{2}\left( u\right) +e_{3}e_{3}\left(
u\right) \\
&&-2ve_{2}\left( u\right) -\frac{1}{2}e_{1}\left( \log \lambda \right)
e_{1}\left( u\right) -\frac{1}{2}e_{3}\left( \log \lambda \right)
e_{3}\left( u\right) .
\end{eqnarray*}%
Moreover, in view of $\left( 5.2\right) $, we get%
\begin{equation*}
e_{1}e_{1}\left( u\right) =e_{1}e_{3}\left( v\right) ,\text{ }%
e_{3}e_{3}\left( u\right) =-e_{3}e_{1}\left( v\right) ,\text{ }%
e_{2}e_{2}\left( u\right) =2ve_{2}\left( u\right) +2ue_{2}\left( v\right) .
\end{equation*}%
Combining $\left( 5.2\right) $, $\left( 5.3\right) $ and the above mentioned
relations, we obtain $\Delta u=0$. Similarly we verify that $v$ is harmonic.
\end{proof}

In the sequel, based on ideas developed by Ramanathan \cite{Ramanathan}, we
intend to give a careful analysis of the set of totally geodesic points%
\begin{equation*}
\mathcal{A}:=\left\{ x\in M^{3}:S\left( x\right) =0\right\} .
\end{equation*}%
Since $f$ is a minimal immersion, it follows that $f$ is real analytic.
Hence, the square $S$ of the length of the second fundamental form is a real
analytic function and the set $\mathcal{A}$ is real analytic. According to
Lojasewicz's Structure Theorem \cite[Theorem 6.3.3, p. 168]{Krantz}, $%
\mathcal{A}$ can be locally decomposed as%
\begin{equation*}
\mathcal{A=V}^{0}\cup \mathcal{V}^{1}\cup \mathcal{V}^{2}\cup \mathcal{V}^{3}
\end{equation*}%
where $\mathcal{V}^{m}$, $0\leq m\leq 3$, is either empty or a finite,
disjoint union of $m$-dimensional subvarieties. A point $x_{0}\in \mathcal{A}
$ is called \textit{smooth point} \textit{of dimension} $m$, if there is a
neighborhood $\Omega $ of $x_{0}$ such that $\Omega \cap \mathcal{A}$ is an $%
m$-dimensional real analytic submanifold of $\Omega $. Otherwise, the point $%
x_{0}$ is said \textit{singular}.

The distribution%
\begin{equation*}
p\in M^{n}\rightarrow \mathcal{D}\left( p\right) :=\ker A\left( p\right) ,
\end{equation*}%
is called the \textit{nullity distribution} of $f$. From our assumptions it
follows that if $x$ is a point of $M^{n}$ then, either $\dim \mathcal{D}%
\left( x\right) =1$ or $\dim \mathcal{D}\left( x\right) =3$. Obviously, the
dimension of the nullity distribution is $3$ precisely at the totally
geodesic points of $M^{3}$.

\begin{lemma}
Let $\Omega $ be an open subset of $M^{3}$ where the dimension of the
nullity distribution is constant. Then $\mathcal{D}$ is differentiable and
involutive. The leaves of $\mathcal{D}$ are totally geodesic in $M^{3}$ and
their images through $f$ are totally geodesic submanifolds in $\mathbb{Q}%
_{c}^{4}$. Moreover, the Gauss map of $f$ is constant along the leaves of $%
\mathcal{D}$.
\end{lemma}

\begin{proof}
Assume at first that the dimension of the nullity distribution in $\Omega $
is $3$. Then all the points of $\Omega $ are totally geodesic and
consequently, the image of $\Omega $ through $f$ lies in a totally geodesic
hypersurface of $\mathbb{Q}_{c}^{4}$. Suppose now that the dimension of $%
\mathcal{D}$ in $\Omega $ is $1$. Let $\gamma :\left( -\varepsilon
,\varepsilon \right) \rightarrow \Omega $ be an integral curve of the vector
field $e_{2}$. From relations $\left( 5.1\right) $ it follows that 
\begin{equation*}
\nabla _{e_{2}}e_{2}=\omega _{21}\left( e_{2}\right) e_{1}+\omega
_{23}\left( e_{2}\right) e_{3}=0.
\end{equation*}%
Thus, $\gamma $ is geodesic in $M^{3}$. Furthermore, from the Gauss formula
we obtain%
\begin{equation*}
\nabla _{e_{2}}df\left( e_{2}\right) =df\left( \nabla _{e_{2}}e_{2}\right)
+\left\langle Ae_{2},e_{2}\right\rangle \xi =0,
\end{equation*}%
and thus, the curve $f\circ \gamma :\left( -\varepsilon ,\varepsilon \right)
\rightarrow \mathbb{Q}_{c}^{4}$ is a geodesic. Finally, from the Weingarten
formula $d\xi =-df\circ A$, we deduce that the Gauss map of $f$ is constant
along the leaves of $\mathcal{D}$.
\end{proof}

The following important result is due of D. Ferus \cite[Lemma 2, p. 311]%
{Ferus}.

\begin{lemma}
Let $\gamma :\left[ 0,b\right] \rightarrow M^{3}$ be a geodesic curve such
that $\gamma \prime \left( t\right) \in \mathcal{D}$, for any $t\in \left[
0,b\right) $, and $\dim \mathcal{D}\left( \gamma \left( t\right) \right) =m$%
, for any $t\in \lbrack 0,b)$. Then, $\dim \mathcal{D}\left( \gamma \left(
b\right) \right) =m$.
\end{lemma}

\begin{proof}
Set $E\left( t\right) =\gamma \prime \left( t\right) =e_{2}\circ \gamma
\left( t\right) $, $t\in \left[ 0,b\right) $, and consider a vector field $Y$
along $\gamma :[0,b)\rightarrow M^{3}$ such that%
\begin{equation}
\left\{ 
\begin{array}{l}
Y\left( 0\right) \perp E\left( 0\right) , \\ 
\\ 
\lbrack Y,E]=0.%
\end{array}%
\right.
\end{equation}%
The existence and uniqueness of $Y$ are immediate since this is a Cauchy
problem for an ordinary linear equation of first order. Note that from $%
\left( 5.4\right) $, we obtain%
\begin{equation}
\frac{dY}{dt}=\nabla _{E}Y=\nabla _{Y}E.
\end{equation}%
Furthermore, we claim that $Y$ is normal along $\gamma :[0,b)\rightarrow
M^{3}$. Indeed,%
\begin{equation*}
\frac{d}{dt}\left\langle E,Y\right\rangle =\left\langle \nabla
_{E}E,Y\right\rangle +\left\langle E,\nabla _{Y}E\right\rangle =\left\langle
E,\nabla _{Y}E\right\rangle =\frac{1}{2}Y\left\langle E,E\right\rangle =0.
\end{equation*}%
Because $\left\langle E\left( 0\right) ,Y\left( 0\right) \right\rangle =0$,
it follows that $\left\langle E\left( t\right) ,Y\left( t\right)
\right\rangle =0$, for any $t\in \lbrack 0,b)$. Differentiating the relation 
$\left( 5.5\right) $, we deduce that%
\begin{equation*}
\frac{d^{2}Y}{dt^{2}}=\nabla _{E}\nabla _{Y}E.
\end{equation*}%
On the other hand we have%
\begin{equation}
R\left( E,Y\right) E=\nabla _{E}\nabla _{Y}E-\nabla _{Y}\nabla _{E}E-\nabla
_{\left[ E,Y\right] }E=\nabla _{E}\nabla _{Y}E=\frac{d^{2}Y}{dt^{2}},
\end{equation}%
where $R$ stands for the Riemannian curvature tensor of $M^{3}$. Moreover,
the Gauss equation implies that%
\begin{equation}
R\left( E,Y\right) E=c\left\{ \left\langle Y,E\right\rangle E-Y\right\} =-cY.
\end{equation}%
Thus, combining the equations $\left( 5.6\right) $ and $\left( 5.7\right) $,
we obtain%
\begin{equation*}
\frac{d^{2}Y}{dt^{2}}+cY=0.
\end{equation*}%
This is a second order ordinary linear\ differential equation with constant
coefficients in $[0,b)$ and thus, the solution can be extended for $t=b$.
Consider now a parallel vector field $Z$ along $\gamma :\left[ 0,b\right]
\rightarrow M^{3}$ and let $G:\left[ 0,b\right] \rightarrow \mathbb{R}$ the
function given by%
\begin{equation*}
G\left( t\right) :=\left\langle AY\left( t\right) ,Z\left( t\right)
\right\rangle ,\quad t\in \left[ 0,b\right] .
\end{equation*}%
Differentiating, we obtain%
\begin{eqnarray*}
\frac{dG}{dt} &=&E\left\langle AY,Z\right\rangle =\left\langle \nabla
_{E}AY,Z\right\rangle =\left\langle \left( \nabla _{E}A\right) Y+A\nabla
_{E}Y,Z\right\rangle \\
&=&\left\langle \left( \nabla _{Y}A\right) E+A\nabla _{E}Y,Z\right\rangle
=\left\langle -A\nabla _{Y}E+A\nabla _{E}Y,Z\right\rangle \\
&=&\left\langle A\left[ E,Y\right] ,Z\right\rangle \\
&=&0.
\end{eqnarray*}%
Thus, $G$ is constant along $\gamma $. Taking $Z\left( 0\right) =Y\left(
0\right) $, we get that%
\begin{equation*}
\left\langle AY,Z\right\rangle \left( \gamma \left( b\right) \right)
=\left\langle AY,Z\right\rangle \left( \gamma \left( 0\right) \right) .
\end{equation*}%
Therefore,%
\begin{equation*}
\dim \mathcal{D}\left( \gamma \left( 0\right) \right) =\dim \mathcal{D}%
\left( \gamma \left( b\right) \right) ,
\end{equation*}%
and the proof is completed.
\end{proof}

Denote by $\mathcal{F}$ the foliation of the nullity distribution in the
open subset $M^{3}-\mathcal{A}$. According to Lemma 5.2, the foliation
consists of geodesics of $M^{3}$. Our aim is to prove that $\mathcal{F}$ can
be analytically extended over the set $\mathcal{A}$. In the following we
shall prove that $\mathcal{A}$ contains only smooth points of dimension $1$.

\begin{lemma}
There are no smooth points of dimension $3$ in $\mathcal{A}$, unless $f$ is
totally geodesic.
\end{lemma}

\begin{proof}
Suppose that $x_{0}$ is a smooth point of $\mathcal{A}$ with dimension $3$.
Then, there exists an open neighborhood $\Omega $ of $M^{3}$ such that $%
x_{0}\in \Omega \subset \mathcal{A}$. From Lemma 5.2, it follows that $%
f:\Omega \rightarrow \mathbb{Q}_{c}^{4}$ is totally geodesic and hence,
according to the Unique Continuation Property, $f:M^{3}\rightarrow \mathbb{Q}%
_{c}^{4}$ must be totally geodesic.
\end{proof}

\begin{lemma}
There are no smooth points of dimension $2$ in $\mathcal{A}$.
\end{lemma}

\begin{proof}
Suppose in contrary that $x_{0}$ is a smooth point of $\mathcal{A}$ with
dimension $2$. Then there exists an open neighborhood $\Omega $ of $M^{3}$
around $x_{0}$ such that $L^{2}=\Omega \cap \mathcal{A}$ is a $2$%
-dimensional imbedded surface. Consider the function $h:M^{3}\rightarrow 
\mathbb{R}$, given by%
\begin{equation*}
h:=\left\{ 
\begin{array}{lll}
\left\langle \xi ,\alpha \right\rangle , & \text{if} & c=0,\smallskip \\ 
\left\langle \xi ,f\left( x_{0}\right) \right\rangle , & \text{if} & c\in
\left\{ -1,1\right\} ,%
\end{array}%
\right.
\end{equation*}%
where $\alpha $ is a vector of $\mathbb{R}^{4}$ such that $h\left(
x_{0}\right) =0$. A simple calculation implies,%
\begin{equation*}
\Delta h+Sh=0,
\end{equation*}%
where $\Delta $ is the Laplacian of $M^{3}$. On the other hand, for any $%
x\in L^{2}$ and $X\in T_{x}M^{3}$ it holds%
\begin{equation*}
dh\left( X\right) =0\text{ and }h\left( x\right) =0.
\end{equation*}%
We claim that $h$ is identically zero. Indeed, let $\left( z;t\right) \in
L^{2}\times \left( -\varepsilon ,\varepsilon \right) $ real analytic
coordinates centered at $x_{0}$. Then, $h$ has a representation of the form,%
\begin{equation*}
h\left( z;t\right) =t^{n}h_{0}\left( z;t\right) ,
\end{equation*}%
where $h_{0}:\Omega \rightarrow \mathbb{R}$, is a non vanishing real
analytic function. Let $z_{1}\in L^{2}$ be a point where $h_{0}\left(
z_{1};0\right) \neq 0$. By continuity there is an $r>0$, such that $h_{0}$
is not zero in the set%
\begin{equation*}
V=B_{\varepsilon }\left( z_{1};0\right) \cap \left( L^{2}\times \left(
0,\varepsilon \right) \right) \text{.}
\end{equation*}%
where $B_{\varepsilon }\left( z_{1};0\right) \subset \Omega $ is the open
ball of radius $r$ centered at $\left( z_{1},0\right) $. Thus, the
restriction of the function $h$ at $V$ has a fixed sign. Assume at first
that $h$ is positive in $V$. Then $\Delta h\leq 0$ on $V$. Also note that $h$
attains a relative minimum at the boundary point $\left( z_{1};0\right) $.
Consequently, by the Hopf Boundary Point Lemma the outward normal derivative
of $h$ at $\left( z_{1};0\right) $ must be non-zero. This contradicts with
the fact that the differential of $h$ vanishes identically at $\left(
z_{1},0\right) $. The same contradiction occurs if we assume that $h$ is
negative in $V$. Thus, $h$ must be identically zero on $\Omega $, which
implies that $f:\Omega \rightarrow \mathbb{Q}_{c}^{4}$ is totally geodesic.
Hence, all the points of $\Omega $ are smooth with dimension 3 which
contradicts to our assumption. Therefore, there are no smooth points of $%
\mathcal{A}$ with dimension 2.
\end{proof}

\begin{lemma}
The set $\mathcal{A}$ does not contain any isolated point.
\end{lemma}

\begin{proof}
Suppose in contrary that $x_{0}$ is an isolated point of $\mathcal{A}$.
Then, there exists an open neighborhood $\Omega $ of $M^{3}$ around $x_{0}$,
such that%
\begin{equation*}
\dim \mathcal{D}\left( x\right) =\left\{ 
\begin{array}{lll}
1, & \text{if} & x\in \Omega -\left\{ x_{0}\right\} , \\ 
3, & \text{if} & x=x_{0}.%
\end{array}%
\right.
\end{equation*}%
Let $\left\{ x_{n}\right\} $ be a sequence of points of $\Omega -\left\{
x_{0}\right\} $ converging to $x_{0}$ and $E_{n}=e_{2}\left( x_{n}\right)
\in T_{x_{n}}M^{3}$ the sequence of unit vectors contained in the nullity
distribution of $f$. By passing to a subsequence, if necessary, we may
assume that there is a vector $E_{0}\in T_{x_{0}}M^{3}$ such that $%
\lim\limits_{n\rightarrow \infty }E_{n}=E_{0}$. Consider now the geodesic
flow $\Phi :\left( -\varepsilon ,\varepsilon \right) \times \Sigma
\rightarrow \Omega $, where%
\begin{equation*}
\Sigma :=\left\{ \left( x,E\left( x\right) \right) \in T_{x}\Omega
:\left\vert E\left( x\right) \right\vert =1\right\} \text{.}
\end{equation*}%
For a fixed pair $\left( x,e_{2}\left( x\right) \right) $, the curve $\gamma
:\left( -\varepsilon ,\varepsilon \right) \rightarrow \Omega $ given by%
\begin{equation*}
\gamma \left( t\right) :=\Phi \left( t,x,e_{2}\left( x\right) \right) ,
\end{equation*}%
is the geodesic passing through $x$ with tangent vector $e_{2}\in T_{x}M^{3}$%
. By continuity, we deduce that%
\begin{equation*}
\lim_{n\rightarrow \infty }\Phi \left( t,x_{n},E_{n}\right) =\Phi \left(
t,x_{0},E_{0}\right) \text{ }\ \text{and \ }\lim_{n\rightarrow \infty }\frac{%
\partial \Phi }{\partial t}\left( t,x_{n},E_{n}\right) =\frac{\partial \Phi 
}{\partial t}\left( t,x_{0},E_{0}\right) ,
\end{equation*}%
where now $\gamma _{0}\left( t\right) :=\Phi \left( t,x_{0},E_{0}\right) $
is the geodesic passing through $x_{0}$ with tangent vector $E_{0}\in
T_{x_{0}}M^{3}$. Moreover, since $\left\Vert A\gamma _{n}^{\prime }\left(
t\right) \right\Vert =0$, again by continuity, we obtain that $\left\Vert
A\gamma _{0}^{\prime }\left( t\right) \right\Vert =0$, for any $t\in \left(
-\varepsilon ,\varepsilon \right) $. So, $\gamma _{n}^{\prime }\left(
t\right) \in \mathcal{D}$ for any $t\in \left( -\varepsilon ,\varepsilon
\right) $. Observe now that%
\begin{equation*}
\dim \mathcal{D}\left( \gamma \left( t\right) \right) =\left\{ 
\begin{array}{ll}
1, & t\neq 0, \\ 
3, & t=0.%
\end{array}%
\right.
\end{equation*}%
This contradicts Lemma 5.3, thus such a point $x_{0}$ does not exist.
\end{proof}

\begin{lemma}
The foliation $\mathcal{F}$ of the nullity distribution can be extended
continuously over the smooth points of $\mathcal{A}$.
\end{lemma}

\begin{proof}
Let $\Omega $ be an open subset of $M^{3}$ such that $\Omega \cap \mathcal{A}
$ is a smooth embedded curve. Consider, as in Lemma 5.6, a sequence of
points $\left\{ x_{n}\right\} \in \Omega -\mathcal{A}$ and a sequence of
unit vectors $E_{n}\in \mathcal{D}$ such that $\lim_{n\rightarrow \infty
}x_{n}=x_{0}\in \mathcal{A}$ and $\lim_{n\rightarrow \infty }E_{n}=E_{0}\in
T_{x_{0}}M^{3}$. Using again the geodesic flow $\Phi $, we deduce\ that the
tangent vector $\gamma _{0}^{\prime }\left( t\right) $ of the geodesic curve%
\begin{equation*}
\gamma _{0}\left( t\right) :=\lim\limits_{n\rightarrow \infty }\Phi \left(
t,x_{n},E_{n}\right) ,
\end{equation*}%
belongs to $\mathcal{D}$, for any $t\in \left( -\varepsilon ,\varepsilon
\right) $. By virtue of Lemma 5.3, $\gamma _{0}$ must lie entirely in $%
\mathcal{A}$. If we consider another sequence $\left\{ \widetilde{x}%
_{n}\right\} $ tending at $x_{0}$ and another sequence of vectors $%
\widetilde{E}_{n}$ tending at $\widetilde{E}_{0}$, then repeating the same
argument as above we deduce that the limit line must lie in $\mathcal{A}$
and so $E_{0}=\pm \widetilde{E}_{0}$. Hence, there is a continuous extension
of the line field $\mathcal{F}$ to the set $\mathcal{A}$.
\end{proof}

\begin{lemma}
The foliation $\mathcal{F}$ of the nullity distribution can be extended
analytically over the smooth points of $\mathcal{A}$.
\end{lemma}

\begin{proof}
Let $TM^{3}$ be\ the tangent bundle of the Riemannian manifold $M^{3}$ and
endow it with the Sasaki metric $\left\langle \text{ },\text{ }\right\rangle
_{\sim }$. In fact, if $X$ and $Y$ are tangent vectors at the time $t=0$ to
the curves $\widetilde{\alpha }\left( t\right) =\left( \alpha \left(
t\right) ,V\left( t\right) \right) $, $\widetilde{\beta }\left( t\right)
=\left( \beta \left( t\right) ,W\left( t\right) \right) $ on $TM^{3}$, then%
\begin{equation*}
\left\langle X,Y\right\rangle _{\sim }:=\left\langle \alpha ^{\prime }\left(
0\right) ,\beta ^{\prime }\left( 0\right) \right\rangle +\left\langle
V_{t}\left( 0\right) ,W_{t}\left( 0\right) \right\rangle ,
\end{equation*}%
where $V_{t}$ and $W_{t}$ stand for the covariant derivatives of the vector
fields $V$, $W$ along the curves $\alpha \left( t\right) $, $\beta \left(
t\right) $ on $M^{3}$. Denote by $T_{1}M^{3}$ the unit tangent bundle of $%
TM^{3}$. Here $T_{1}M^{3}$ is regarded as a submanifold of $TM^{3}$ equipped
with the induced Riemannian metric by the Sasaki metric of $TM^{3}$. Assume
now that $\Omega $ is an open subset of $M^{3}-\mathcal{A}$ and that$\
\left\{ e_{1},e_{2},e_{3}\right\} $ is\ a local frame of principal
directions on $\Omega $, such that%
\begin{equation*}
Ae_{1}=\lambda e_{1}\text{, \ }Ae_{2}=0\text{, \ }Ae_{3}=-\lambda e_{3}.
\end{equation*}%
Define the smooth map $F:\Omega \rightarrow T_{1}M^{3}$, by%
\begin{equation*}
F\left( x\right) :=\left( x,e_{2}\left( x\right) \right) ,\quad x\in \Omega .
\end{equation*}%
We claim that the map $F$ is harmonic. Indeed, from the relations $\left(
5.1\right) $ and $\left( 5.2\right) $ we can compute the rough Laplacian $%
\widetilde{\Delta }F$ of $F$,%
\begin{eqnarray*}
\widetilde{\Delta }F &=&\dsum\limits_{i=1}^{3}\left\{ \nabla _{e_{i}}\nabla
_{e_{i}}F-\nabla _{\nabla _{e_{i}}e_{i}}F\right\} \\
&=&\nabla _{e_{1}}\nabla _{e_{1}}e_{2}+\nabla _{e_{3}}\nabla
_{e_{3}}e_{2}-\nabla _{\nabla _{e_{1}}e_{1}}e_{2}-\nabla _{\nabla
_{e_{1}}e_{1}}e_{2} \\
&=&-2\left( u^{2}+v^{2}\right) F \\
&=&-2\left( \left\Vert \nabla _{e_{1}}F\right\Vert ^{2}+\left\Vert \nabla
_{e_{3}}F\right\Vert ^{2}\right) F.
\end{eqnarray*}%
Hence, the map $F$ satisfies the differential equation%
\begin{equation*}
\widetilde{\Delta }F+2\left\Vert \nabla F\right\Vert ^{2}F=0,
\end{equation*}%
which is precisely the Euler-Lagrange equation for the energy of the map $F$
(cf., \cite[Proposition 1.1]{Wood}). Consequently, $F$ is harmonic.
Moreover, according to Lemma 5.7, the map$\ F$ can be extended continuously
on $\mathcal{A}$. Consequently, from Theorem 2.2, we deduce that $F$ can be
extended analytically over $\mathcal{A}$. Hence, the foliation $\mathcal{F}$
of the nullity distribution can be extended analytically over the smooth
points of $\mathcal{A}$.
\end{proof}

\begin{lemma}
There are no singular points on the set $\mathcal{A}$.
\end{lemma}

\begin{proof}
Assume in contrary that $x_{0}$ is a singular point of $\mathcal{A}$. From
Lemmata 5.4, 5.5 and 5.6, it follows that $\mathcal{A}$ contains algebraic
curves. The singular points of an algebraic curve are either cusps,
intersection points or combination of those. Furthermore, the singular
points are isolated. Suppose at first that $x_{0}$ is an intersection point.
Let $\Omega $ be an open neighborhood around $x_{0}$, such that the
restriction of $f$ on $\Omega $ is injective, and consider a transversal $%
L^{2}$ to the foliation $\mathcal{F}$ on $U-\left\{ x_{0}\right\} $. Since
the images of the integral curves of $\mathcal{D}$ via $f$ are geodesics
lines in $\mathbb{Q}_{c}^{4}$, the map $f$ can be represented by the formula,%
\begin{equation*}
f\left( z,t\right) =\left\{ 
\begin{array}{ll}
\cos tf\left( z,0\right) +\sin tdf\left( e_{2}\right) \left( z,0\right) , & 
c=\ \ 1,\smallskip \\ 
f\left( z,0\right) +tdf\left( e_{2}\right) \left( z,0\right) , & c=\ \
0,\smallskip \\ 
\cosh tf\left( z,0\right) +\sinh tdf\left( e_{2}\right) \left( z,0\right) ,
& c=-1.%
\end{array}%
\right.
\end{equation*}%
where $z\in L^{2}$. From the above description, it follows that there are
points $\left( z_{1},t_{1}\right) $ and $\left( z_{2},t_{2}\right) $ such
that $f\left( z_{1},t_{1}\right) =f\left( x_{0}\right) =f\left(
z_{2},t_{2}\right) $, which is a contradiction. Suppose now that $x_{0}$ is
a cusp. Since the images of the leaves of the nullity foliation via $f$ are
geodesics in $\mathbb{Q}_{c}^{4}$, we deduce that there must be a point $%
x_{1}\in \Omega $, $x_{1}\neq x_{0}$, where the leaves of $\mathcal{F}$
intersect. This leads again to a contradiction. Thus, there are no singular
points on $\mathcal{A}$.
\end{proof}

\section{Proofs of theorems}

\begin{proof}[Proof of Theorem 1.1]
Assume that $f$ is not totally geodesic. It is enough to examine only the
case where the set $\mathcal{A}$ of the totally geodesic points of $f$ is
non empty. By virtue of the results in Section 5, the set$\ \mathcal{A}$
contains geodesic lines and the foliation $\mathcal{F}$ of the nullity
distribution $\mathcal{D}$ is well defined and regular in the sense of
Palais \cite{Palais}. Consider a point $x_{0}\in \mathcal{A}$ and a
coordinate system $\left( x_{1},x_{2},x_{3}\right) $ on $\Omega \subset
M^{3} $, around $x_{0}$, such that $\partial /\partial x_{2}=e_{2}$. Denote
by $M^{2}=\Omega /\mathcal{F}$ the space of leaves of $\mathcal{F}$ on $%
\Omega $ and by $\pi :\Omega \rightarrow M^{2}$ the quotient projection. It
is well known that $M^{2}$ can be equipped with a differentiable structure
making $\pi $ a submersion. Endow $M^{2}$ with the induced conformal
structure by the Riemannian metric $\left\langle \text{ },\text{ }%
\right\rangle _{0}$ which makes $\pi $ a Riemannian submersion.\medskip

\noindent \textit{Euclidean Case:}\textbf{\ }Since\ the unit normal $\xi $
of the hypersurface remains constant along each leaf of $\mathcal{F}$, we
may define a smooth map $g:M^{2}\rightarrow \mathbb{S}^{3}$ so that $g\circ
\pi :=\xi $. Because $\xi $ is harmonic and $\pi $ is a Riemannian
submersion, it follows that $g$ is also harmonic. Consider a smooth
transversal $L^{2}$ to the leaves of $\mathcal{F}$ passing through the point 
$x_{0}$. Using the Weingarten formula, we deduce that the induced by $g$
metric tensor $\left\langle \text{ },\text{ }\right\rangle _{g}$ on $M^{2}$
is given by%
\begin{equation*}
\left\langle \text{ },\text{ }\right\rangle _{g}\circ \pi \left\vert
_{L}\right. =\lambda ^{2}\left\vert _{L}\right. \left\langle \text{ },\text{ 
}\right\rangle _{0}\circ \pi \left\vert _{L}\right. ,
\end{equation*}%
where $\pi \left\vert _{L}\right. $ and $\lambda ^{2}\left\vert _{L}\right. $
is the restriction of $\pi $ and$\ \lambda ^{2}$ to the transversal $L^{2}$,
respectively. Consequently, $\pi \left( x_{0}\right) $ is a branch point of $%
g$ and thus $g$ is a branched minimal surface. Moreover, the normal bundle
of $g$ is spanned by the vector field $\eta $ given by $\eta \circ \pi
:=df\left( e_{2}\right) $. Denote by $\sigma :\Omega \rightarrow \mathbb{R}$%
, $\sigma :=\left\langle \xi ,f\right\rangle $, the support function of $f$.
Because $e_{2}\left( \sigma \right) =0$, it follows that $\sigma $ is
constant along the leaves of the foliation $\mathcal{F}$. Thus, the function 
$\gamma :M^{2}\rightarrow \mathbb{R}$ given by $\gamma \circ \pi :=$ $\sigma 
$ is well defined and smooth. We may choose the transversal $L^{2}$ in such
a way that the frame $\left\{ E_{1}:=e_{1}\left( x\right) ,\text{ }%
E_{3}:=e_{3}\left( x\right) \right\} $ spans $T_{x}L^{2}$, at a point $x\neq
x_{0}$. Then, $\left\{ X_{1}:=\frac{1}{\lambda \left( x\right) }d\pi \left(
E_{1}\right) \text{, }X_{3}:=\frac{1}{\lambda \left( x\right) }d\pi \left(
E_{3}\right) \right\} $ constitute an orthonormal base at $T_{\pi \left(
x\right) }M^{2}$ with respect to $\left\langle \text{ },\text{ }%
\right\rangle _{g}$ and%
\begin{eqnarray*}
dg\left( \nabla \gamma \right) \circ \pi \left( x\right) &=&dg\left\{ \left(
X_{1}\gamma \right) X_{1}+\left( X_{2}\gamma \right) X_{2}\right\} \\
&=&\left\langle f\left( x\right) ,df\left( E_{1}\right) \right\rangle
df\left( E_{1}\right) +\left\langle f\left( x\right) ,df\left( E_{3}\right)
\right\rangle df\left( E_{3}\right) \\
&=&f^{\top }\left( x\right) ,
\end{eqnarray*}%
where $f^{\top }$ stands for the tangential component of $f$ on $df\left(
T_{x}M^{3}\right) $. However, the right hand side of the above equality is
well defined for any $x\in \Omega $ which implies that $dg\left( \nabla
\gamma \right) $ is well defined even on the branch point. Analyzing $f$
with respect to the orthonormal basis $\left\{ df\left( e_{1}\right)
,df\left( e_{2}\right) ,df\left( e_{3}\right) ,\xi \right\} $, we deduce that%
\begin{equation*}
f=\gamma \left( g\circ \pi \right) +dg\left( \nabla \gamma \right) \circ \pi
+\left\langle f,df\left( e_{2}\right) \right\rangle df\left( e_{2}\right) .
\end{equation*}%
Moreover, $f=\Psi _{\gamma ,g}\circ T$ on $\Omega $, where $T:\Omega
\rightarrow \mathcal{N}\left( g\right) $ is the map defined by%
\begin{equation*}
T\left( x\right) =\left( \pi \left( x\right) ,\left\langle f\left( x\right)
,df\left( e_{2}\left( x\right) \right) \right\rangle df\left( e_{2}\left(
x\right) \right) \right) ,\text{ \ }x\in \Omega .
\end{equation*}%
\noindent \textit{Spherical Case:} As in the Euclidean case, the map $%
g:M^{2}\rightarrow \mathbb{S}^{4}$ given by $g\circ \pi :=\xi $ is well
defined and harmonic. The induced metric tensor by $g$ $\left\langle \text{ }%
,\text{ }\right\rangle _{g}$ on $M^{2}$ is given by%
\begin{equation*}
\left\langle \text{ },\text{ }\right\rangle _{g}\circ \pi \left\vert
_{L}\right. =\lambda ^{2}\left\vert _{L}\right. \left\langle \text{ },\text{ 
}\right\rangle _{0}\circ \pi \left\vert _{L}\right. ,
\end{equation*}%
where $L^{2}$ a smooth transversal to the leaves of $\mathcal{F}$ passing
through the point $x_{0}$ and $\lambda ^{2}\left\vert _{L}\right. $ is the
restriction of $\lambda ^{2}$ to $L^{2}$. Hence, $g:M^{2}\rightarrow \mathbb{%
S}^{4}$ is a branched minimal surface. Moreover, the normal bundle of $g$ is
spanned by the vector fields $\eta _{1}$, $\eta _{2}$ given by $\eta
_{1}\circ \pi :=df\left( e_{2}\right) $ and $\eta _{2}\circ \pi :=f$.
Furthermore, $f=\Psi _{g}\circ T$ on $\Omega $, where the map $T:\Omega
\rightarrow \mathcal{N}_{1}\left( g\right) $ is defined by $T\left( x\right)
=\left( \pi \left( x\right) ,f\left( x\right) \right) $, $x\in \Omega $%
.\medskip

\noindent \textit{Hyperbolic Case: }The proof is the same as in the
spherical case.
\end{proof}

\begin{proof}[Proof of Theorem 1.2]
The assumption in \cite{Hasanis1, Hasanis2, Hasanis3} that $f$ is free of
totally geodesic points was used only to show that the local harmonic
functions%
\begin{equation*}
u=\left\langle \nabla _{e_{3}}e_{1},e_{2}\right\rangle \text{ \ and \ }%
v=\left\langle \nabla _{e_{1}}e_{1},e_{2}\right\rangle ,
\end{equation*}%
can be defined globally on $M^{3}$. However, due to results in Section 5,
the local functions $u$ and $v$ can be extended analytically on $M^{3}$ even
in the case when $f$ allows totally geodesic points. Indeed, suppose that $f$
is not totally geodesic. Then, the foliation $\mathcal{F}$ of the nullity
distribution $\mathcal{D}$ can be extended analytically on the set $\mathcal{%
A}$ of the totally geodesic points of $f$. Without loss of generality we may
assume that there is a global unit section $e_{2}\in \mathcal{D}$, since
otherwise we can pass to the $2$-fold covering space%
\begin{equation*}
\widetilde{M}^{3}:=\left\{ \left( x,\varepsilon _{x}\right) :x\in
M^{3},\varepsilon _{x}\in \mathcal{D}\left( x\right) \right\} .
\end{equation*}%
Since the section $e_{2}$ is globally defined, the local functions $u$\ and $%
v$ can be extended analytically to harmonic functions defined on all of $%
M^{3}$. Consequently, following the same arguments as in \cite{Hasanis1,
Hasanis2, Hasanis3}, we can complete the proof.
\end{proof}

\section{Acknowledgements}

The author is grateful to professor J.-H. Eschenburg for his helpful
comments and suggestions. Moreover, the author would like also to thank
professors A. Tertika and S. Filippa for many fruitful discussions.

\end{document}